\documentclass{article}
\usepackage{graphicx} 

\usepackage{amsmath,amssymb,amsthm,mathtools,fullpage}

\newtheorem{theorem}{Theorem}
\newtheorem{lemma}[theorem]{Lemma}
\newtheorem{corollary}[theorem]{Corollary}
\theoremstyle{remark}
\newtheorem{remark}[theorem]{Remark}
\usepackage{titlesec}
\titleformat*{\section}{\large\bfseries}
\titleformat*{\subsection}{\normalsize\bfseries}
\usepackage{hyperref}

\title{\textbf{Singular Value Counts and Metric Entropy\\
in Calder\'on's Inverse Conductivity Problem}}

\author{%
  Leonid Sarieddine\textsuperscript{1}
  \quad and \quad
  Raul Jimenez\textsuperscript{1,2}\\[0.6em]
  {\small \textsuperscript{1}ICCUB, University of Barcelona,
  Mart\'i i Franqu\`es 1, 08028 Barcelona, Spain}\\
  {\small \textsuperscript{2}ICREA, Passeig Llu\'is Companys 23,
  08010 Barcelona, Spain}
}

\date{}
\date{}

\begin{document}
\maketitle

\begin{abstract}

We ask how much of the conductivity in the interior of a body can be
resolved from boundary measurements at a given noise level. Let $\Omega$
be a bounded Lipschitz domain in $\mathbb R^d$, with $d \geq 2$, and
let $E$ be a nonempty open set compactly contained in $\Omega$. We
measure changes in the Dirichlet-to-Neumann map in the Hilbert--Schmidt
norm on bilinear forms on $H^{1/2}(\partial\Omega)$. For the map
linearized at conductivity $1$, acting on perturbations in $L^2(E)$,
we prove that the number of singular values above $\varepsilon$ is of
order $(\log(1/\varepsilon))^d$ as $\varepsilon$ tends to zero.

For the full non-linear problem we fix a background Lipschitz conductivity $\beta$ and consider all conductivities $\gamma$ that
agree with $\beta$ outside $E$ and satisfy fixed ellipticity bounds.
Let $N_\varepsilon$ be the smallest number of Hilbert--Schmidt balls
of radius $\varepsilon$ needed to cover the set of differences
$\Lambda_\gamma-\Lambda_\beta$. We show that, for all sufficiently
small $\varepsilon$, $
c(\log(1/\varepsilon))^{d+1}
\leq \log N_\varepsilon
\leq C(\log(1/\varepsilon))^{2d-1}$.
The same bounds hold in the operator norm from $H^{1/2}(\partial\Omega)$ to $H^{-1/2}(\partial\Omega)$. In two dimensions the bounds agree, giving the sharp order
$(\log(1/\varepsilon))^3$ for the entropy of the nonlinear data.
    
\end{abstract}

\section{Introduction and main results}

Calder\'on's problem \cite{Calderon2006} asks whether the conductivity $\gamma$, inside a domain
$\Omega\subset\mathbb{R}^d$, $d\ge 2$, can be recovered from the Dirichlet-to-Neumann (DtN)
map $\Lambda_\gamma$, which encodes all voltage and current measurements on the boundary denoted by
$\partial\Omega$. It is the mathematical model of electrical impedance tomography
\cite{Borcea2002,Uhlmann2009}. On a simply connected planar domain the full boundary map determines a
bounded, uniformly positive scalar conductivity~\cite{AstalaPaivarinta2006}.
In dimensions $d\geq3$, uniqueness was proved for smooth conductivities
by Sylvester and Uhlmann~\cite{SylvesterUhlmann1987} and for Lipschitz
conductivities by Caro and Rogers~\cite{CaroRogers2016}. For bounded
measurable conductivities in three dimensions, a recent preprint gives
two distinct examples on a ball that equal one near the boundary and
have the same boundary map~\cite{OpenAIConductivity2026}.
Under suitable a priori regularity bounds, the inverse conductivity
map satisfies logarithmic stability estimates~\cite{Alessandrini1988,ClopFaracoRuiz2010}.
See also the survey~\cite{Alessandrini2007}.
Mandache's examples rule out H\"older stability in general, even under
a fixed smoothness bound~\cite{Mandache2001}. His proof compares the
number of separated coefficients with the number of balls needed to
cover their boundary maps. This argument has been developed further
in~\cite{DiCristoRondi2003,KowWang2022,KochRulandSalo2025}. For
piecewise constant conductivities on a prescribed partition, Lipschitz
stability holds~\cite{AlessandriniVessella2005}, although the stability
constant can grow rapidly with the number of
subdomains~\cite{Rondi2006}.

We ask how many different boundary responses can be distinguished at
a given noise level when the conductivity varies only in the interior.
For the linearized map we count singular values. For the nonlinear map
we estimate the number of balls needed to cover its image.

Let $E\Subset\Omega$ be a nonempty open set, let $H=H^{1/2}(\partial\Omega)$, and measure
boundary data in the Hilbert--Schmidt norm on bilinear forms on $H\times H$. Our first
result concerns the linearization of $\gamma\mapsto\Lambda_\gamma$ at the constant
conductivity,
\[
  (Aq)(f,g)=\int_E q\,\nabla u_f\cdot\nabla u_g\,dx ,
\]
where $u_f$ is the harmonic extension of $f$; this is the operator introduced by
Calder\'on \cite{Calderon2006}. We show that $A\colon L^2(E)\to\mathrm{HS}$ is compact and
that the number $N(\tau)$ of eigenvalues of $A^*A$ exceeding $\tau$ satisfies
\[
  N(\tau)\asymp\log(1/\tau)^d \qquad (\tau\downarrow 0),
\]
see Theorem~\ref{linear theorem}. Equivalently, the singular values of $A$ satisfy
$c\,e^{-Cn^{1/d}}\le s_n(A)\le C\,e^{-cn^{1/d}}$. Thus $N(\tau)$ counts the unit $L^2(E)$ singular modes whose linearized boundary response has Hilbert–Schmidt norm at least $\sqrt{\tau}$.

Our second result concerns the full non-linear map. Fix $0<\lambda<L$ and a Lipschitz
background conductivity $\beta$ satisfying $\lambda \le \beta \le L$, and let $K_\beta$ be the set of measurable conductivities
with $\lambda\le\gamma\le L$ that agree with $\beta$ outside $E$; no regularity is assumed
inside $E$. Then the Kolmogorov $\varepsilon$-entropy \cite{KT} of
the set of DtN maps, in the Hilbert--Schmidt norm, satisfies
$$
c\,\log(1/\varepsilon)^{d+1}\ \le\
H_\varepsilon\bigl(\{\Lambda_\gamma-\Lambda_\beta:\gamma\in K_\beta\}\bigr)\ \le\
C\,\log(1/\varepsilon)^{2d-1}
$$
for small $\varepsilon$, see Theorems~\ref{nonlinear lower bound} and~\ref{nonlinear upper bound}; in particular it
is of order $\log(1/\varepsilon)^3$ in the plane. The upper exponent $2d-1$ also appears in Mandache's covering
estimate~\cite{Mandache2001}. For conductivities in a ball that equal
one outside an interior cube, Rondi obtained the same exponent in the
usual operator norm~\cite{Rondi2006}. Our upper bound holds in
the Hilbert--Schmidt norm on bounded Lipschitz domains with a prescribed
exterior conductivity.

Several lines of work are closely related. On the unit disc, Allers and Santosa
\cite{AllersSantosa1991} reduced the linearized problem to a moment problem and concluded
that any stable reconstruction has limited resolution. For the planar Neumann-to-Dirichlet map, the derivative at the constant
conductivity takes $L^2$ perturbations to Hilbert--Schmidt operators on
mean-free $L^2(\partial\Omega)$, even when the perturbations reach the
boundary~\cite{BischHirvensaloHyvonen2025}.
Garde and Hyv\"onen proved Lipschitz stability on certain
infinite-dimensional subspaces~\cite{GardeHyvonen2024}. For the
linearized Schrodinger problem on the disk, Lionheart obtained an
explicit spectral representation using multiplicity-normalized
data~\cite{Lionheart2026}. Lionheart, Pidcock and Sebu treated radial
potentials and observed that the normal operator becomes compact when
the potential is fixed near the boundary~\cite{LionheartPidcockSebu2026}.
Here we obtain two-sided singular-value bounds for the conductivity
derivative on interior-supported $L^2$ perturbations in every
dimension $d\geq2$. For radial potentials the linear part is known to carry
the whole instability of the Calder\'on problem \cite{DaudeMaciaMeronoNicoleau2026}. Limits
of resolution in impedance imaging have been studied through distinguishability
\cite{Isaacson1986,GardeKnudsen2017}, the information content of the data
\cite{Dobson1992,AdlerYoumaranLionheart2008} and depth-dependent stability
\cite{NagayasuUhlmannWang2009,AlessandriniScapin2017}, and the regularized D-bar method
truncates the nonlinear Fourier transform at a radius proportional to
$\log(1/\varepsilon)$ \cite{KnudsenLassasMuellerSiltanen2009}. At positive frequency the
singular values of linearized scattering problems first stay roughly constant on a set of
stably recoverable features and only then decay exponentially \cite{KowSaloZou2025}; see
\cite{LiimatainenSaloTrad2026} for low Fourier modes in the nonlinear problem. Finally,
logarithmic convergence rates are the best possible also in the statistical Calder\'on
problem \cite{AbrahamNickl2019}.

\section{Singular-value counting for the linearized map}
\label{sec:linearized-counting}

Throughout, $\Omega\subset\mathbb R^d$, $d\geq2$, is a bounded connected
Lipschitz domain and $E\Subset\Omega$ is a nonempty open set. We equip
$H=H^{1/2}(\partial\Omega)$ with the complexification of a fixed real
Sobolev inner product. Inner products are linear in the first argument.

For a bounded, uniformly positive conductivity $\gamma$ and $f\in H$,
let $u_f^\gamma\in H^1(\Omega)$ be the weak solution of
\[
\nabla\cdot(\gamma\nabla u_f^\gamma)=0 \quad\text{in }\Omega,
\qquad u_f^\gamma|_{\partial\Omega}=f.
\]
We write $u_f:=u_f^1$ for the harmonic extension of $f$. The
Dirichlet-to-Neumann form is
\[
\Lambda_\gamma(f,g)=\int_\Omega
\gamma\,\nabla u_f^\gamma\cdot\nabla u_g^\gamma\,dx,
\]
where the dot denotes the complex bilinear extension of the Euclidean
dot product. We write $\mathrm{HS}$ for the Hilbert space of
Hilbert--Schmidt bilinear forms on $H\times H$, with norm
$\|B\|_{\mathrm{HS}}^2=\sum_{j,k}|B(e_j,e_k)|^2$ for any orthonormal
basis $(e_j)$ of $H$~\cite{Kostenko}.
We first take a real perturbation $q\in L^\infty(E)$ and extend it by
zero outside $E$, so that $1+tq$ remains uniformly positive for $|t|$
small. Differentiating at $t=0$ gives
$$
(Aq)(f,g)
:=\left.\frac{d}{dt}\right|_{t=0}\Lambda_{1+tq}(f,g)
=\int_E q\,\nabla u_f\cdot\nabla u_g\,dx.
$$
Indeed, the boundary values are fixed, so the derivatives of the
solutions have zero trace. The two terms involving these derivatives
therefore vanish by the weak equations for $u_f$ and $u_g$.

To study the singular values of $A$, we take $L^2(E)$ as the space of
perturbations. This gives a Hilbert space without imposing any smoothness
on $q$. Here $L^2(E)$ denotes the complex space, and $A$ is extended complex-linearly.
Since the harmonic extension of a real function is real, $A$ is the
complexification of its restriction to real-valued $q$; in particular, both
have the same singular values. We first show that $A$, initially defined for bounded
perturbations, extends to a bounded operator from $L^2(E)$ to $\mathrm{HS}$.

Choose an open set $U$ with $\overline E\subset U\Subset\Omega$.
Interior estimates for harmonic functions, together with the
Dirichlet estimate, give
$\|\nabla u_f\|_{C^m(U)}\le C_{m,U,\Omega}\|f\|_H$ for every
$m\ge0$; [see \cite{Evans2010},
\cite{RS}]. Thus, for $x\in U$ and $i=1,\ldots,d$,
the functional $L_{i,x}(f):=\partial_i u_f(x)$ is bounded on $H$.
Define
\[
B_x(f,g)
:=
\nabla u_f(x)\cdot\nabla u_g(x)
=
\sum_{i=1}^d L_{i,x}(f)L_{i,x}(g).
\]
Each $B_x$ is a finite-rank bilinear form, and
$\|B_x\|_{\mathrm{HS}}
\le\sum_{i=1}^d\|L_{i,x}\|_{H^*}^2\le C$.
The second-derivative estimate also gives
$\|L_{i,x}-L_{i,y}\|_{H^*}\le C|x-y|$ whenever the segment
joining $x$ and $y$ lies in a sufficiently small ball in $U$.
Consequently, $x\mapsto B_x$ is continuous in the
Hilbert--Schmidt norm, hence strongly measurable. For every $q\in L^2(E)$, the Bochner integral $\int_E q(x)B_x\,dx$
is well defined in $\mathrm{HS}$, since
\[
\int_E |q(x)|\,\|B_x\|_{\mathrm{HS}}\,dx
\le C|E|^{1/2}\|q\|_{L^2(E)}.
\]
For bounded $q$, evaluating this integral at $(f,g)$ gives
$\int_E q\,\nabla u_f\cdot\nabla u_g\,dx$, which is the formula
for $Aq$ obtained above. The integral therefore defines a bounded
linear extension $A:L^2(E)\to\mathrm{HS}$. This extension is unique
because $L^\infty(E)$ is dense in $L^2(E)$.

We will show below that $A$ is in fact a compact operator and hence by the spectral theorem of compact self adjoint operators $A^*A$ has a countable set of eigenvalues with $0$ being the only possible accumulation point. We define the following counting function

$$N(\tau):= \# \{n: \lambda_n(A^*A) \ge \tau\}$$
where $\lambda_n$ denotes the eigenvalues of $A^*A$. Note that for each positive $\tau$, $N(\tau)$ is a finite number since by compactness of $A$ only finitely many eigenvalues are above $\tau$. $N(\tau)$ can be interpreted as counting the number of $L^2$ singular modes whose boundary response has norm at least $\sqrt{\tau}$.

\begin{theorem}\label{thm:linearized-counting}
As $\tau\downarrow0$,
\[
N(\tau)\asymp\bigl(\log(1/\tau)\bigr)^d.
\]
The constants depend only on $d,\Omega,E$ and the chosen norm
on $H$.
\label{linear theorem}
\end{theorem}

\subsection{Interior analyticity and the upper bound}

We use the following elementary consequence of Taylor approximation.

\begin{lemma}\label{lem:linear-analytic-approximation}
Let $E\Subset U\subset\mathbb R^d$, where $U$ is open, and let
$\mathcal G$ be a Hilbert space. Suppose that
$F:U\to\mathcal G$ is real analytic and satisfies
$\|\partial^\alpha F(x)\|_{\mathcal G}
\le C_0R^{|\alpha|}\alpha!$ for every $x\in U$ and every
multi-index $\alpha$. Define $T:\mathcal G\to L^2(E)$ by
$(Th)(x):=\langle h,F(x)\rangle_{\mathcal G}$.
Then $T$ is compact and $s_n(T)\le Ce^{-cn^{1/d}}$ for all
$n\ge1$.
\end{lemma}

\begin{proof}
Choose $r>0$ small enough that $Rr\le1/4$ and every cube of
half-side length $r$ centered on $\overline E$ is contained in $U$.
On such a cube centered at $x_j$, let
\[
P_{m,j}(x)
:=
\sum_{|\alpha|\le m}
\frac{\partial^\alpha F(x_j)}{\alpha!}(x-x_j)^\alpha.
\]
For $|x-x_j|_\infty\le r$, the Taylor series and the derivative
bounds give
\begin{align}
\|F(x)-P_{m,j}(x)\|_{\mathcal G}
&\le C_0\sum_{|\alpha|>m}(Rr)^{|\alpha|} \notag\\
&=
C_0\sum_{k>m}\binom{k+d-1}{d-1}(Rr)^k
\le C\,2^{-m}.
\label{eq:linear-taylor-error}
\end{align}
Indeed, $\binom{k+d-1}{d-1}\le C_d(k+1)^{d-1}$ and
$Rr\le1/4$, so the polynomial factor is absorbed by the
geometric decay.

Cover $\overline E$ by finitely many such cubes
$Q_1,\ldots,Q_J$, and choose a measurable partition
$E=\bigsqcup_{j=1}^J E_j$ with $E_j\subset Q_j$.
Define $P_m(x):=P_{m,j}(x)$ on $E_j$ and
$(T_mh)(x):=\langle h,P_m(x)\rangle_{\mathcal G}$.
The range of $T_m$ lies in the span of
$\mathbf1_{E_j}(x)(x-x_j)^\alpha$ with $|\alpha|\le m$,
so $\operatorname{rank}T_m\le J\binom{m+d}{d}\le C(m+1)^d$.
Moreover, \eqref{eq:linear-taylor-error} gives
\[
\|(T-T_m)h\|_{L^2(E)}^2
\le
\int_E\|h\|_{\mathcal G}^2
       \|F(x)-P_m(x)\|_{\mathcal G}^2\,dx
\le C^2|E|\,2^{-2m}\|h\|_{\mathcal G}^2.
\]
Thus $T$ is a norm limit of finite-rank operators.

For compact operators between Hilbert spaces,
$s_{n+1}(T)=\inf_{\operatorname{rank}S\le n}\|T-S\|$;
see \cite{Kostenko}. In particular,
\[
s_{J\binom{m+d}{d}+1}(T)
\le \|T-T_m\|
\le C|E|^{1/2}2^{-m}.
\]
Choosing $m$ proportional to $n^{1/d}$, with a sufficiently
small proportionality constant, and using monotonicity of
the singular values proves the claim.
\end{proof}

Now we verify the factorial bounds needed to apply
Lemma~\ref{lem:linear-analytic-approximation} to $A^*$.
Keep the neighborhood $U$ chosen above and fix
$\rho:=\tfrac12\operatorname{dist}(\overline U,\partial\Omega)>0$.
Then $B(x,\rho)\subset\Omega$ for every $x\in U$.

For a multi-index $\alpha$ and a coordinate vector $e_i$,
$\partial_x^\alpha L_{i,x}(f)=\partial^{\alpha+e_i}u_f(x)$.
These derivatives exist in the norm of $H^*$: applying the
interior estimates to one additional derivative controls the
difference quotients uniformly over $\|f\|_H\le1$.

Put $p:=|\alpha|$ and $k:=p+1$. The local harmonic derivative
estimate \cite{Evans2010} gives
\begin{equation}\label{eq:linear-harmonic-derivatives}
|\partial^{\alpha+e_i}u_f(x)|
\le
\frac{C_k}{\rho^{d+k}}\|u_f\|_{L^1(B(x,\rho))},
\qquad
C_k\le C_d(C_dk)^k.
\end{equation}
By Cauchy--Schwarz and the Dirichlet estimate,
\[
\|u_f\|_{L^1(B(x,\rho))}
\le |B(0,\rho)|^{1/2}\|u_f\|_{L^2(\Omega)}
\le C\|u_f\|_{H^1(\Omega)}
\le C'\|f\|_H.
\]
Substituting into \eqref{eq:linear-harmonic-derivatives} and
taking the supremum over $\|f\|_H\le1$, we obtain
$\|\partial_x^\alpha L_{i,x}\|_{H^*}
\le C\rho^{-p-1}C_{p+1}$, where the fixed factor $\rho^{-d}$
has been absorbed into $C$.

Stirling's formula gives $(C_dk)^k\le C_1^k k!$, so after
changing constants,
$\|\partial_x^\alpha L_{i,x}\|_{H^*}\le CR_1^p(p+1)!$.
To pass from the ordinary factorial to the multi-index factorial,
note that $p!\le d^p\alpha!$ by the multinomial theorem and
$p+1\le2^p$. Thus $(p+1)!\le(2d)^p\alpha!$.
Absorbing $(2d)^p$ into $R_1^p$ gives
\begin{equation}\label{eq:linear-functional-factorial}
\|\partial_x^\alpha L_{i,x}\|_{H^*}
\le CR^{|\alpha|}\alpha!,
\qquad x\in U.
\end{equation}

These bounds also establish analyticity in the operator norm.
Indeed, Taylor's formula along a segment in $U$ bounds the
remainder after degree $m$ by
$C\binom{m+d}{d-1}(R|h|_\infty)^{m+1}$ for the increment $h$,
uniformly over $\|f\|_H\le1$. For $R|h|_\infty<1$, this tends
to zero. Hence the Taylor series for $x\mapsto L_{i,x}$
converges in $H^*$.

Since $B_x=\sum_{i=1}^dL_{i,x}\otimes L_{i,x}$, Leibniz's
rule gives
\[
\partial_x^\alpha B_x
=
\sum_{i=1}^d\sum_{\nu\le\alpha}
\binom{\alpha}{\nu}
(\partial_x^\nu L_{i,x})
\otimes
(\partial_x^{\alpha-\nu}L_{i,x}).
\]
Using \eqref{eq:linear-functional-factorial} and the
Hilbert--Schmidt norm of a rank-one form, we obtain
\begin{align}
\|\partial_x^\alpha B_x\|_{\mathrm{HS}}
&\le
CR^{|\alpha|}
\sum_{\nu\le\alpha}
\binom{\alpha}{\nu}\nu!(\alpha-\nu)! \notag\\
&=
CR^{|\alpha|}\alpha!\prod_{j=1}^d(\alpha_j+1)
\le C(2R)^{|\alpha|}\alpha!.
\label{eq:linear-bilinear-factorial}
\end{align}
Thus $x\mapsto B_x$ is real analytic as an
$\mathrm{HS}$-valued map, with uniform factorial bounds on $U$.

Finally, for $q\in L^2(E)$ and $G\in\mathrm{HS}$,
\[
\langle Aq,G\rangle_{\mathrm{HS}}
=
\int_E q(x)\langle B_x,G\rangle_{\mathrm{HS}}\,dx.
\]
Since the inner products are linear in the first argument,
this identifies the adjoint as
$(A^*G)(x)=\langle G,B_x\rangle_{\mathrm{HS}}$.
Lemma~\ref{lem:linear-analytic-approximation}, applied with
$\mathcal G=\mathrm{HS}$ and $F(x)=B_x$, therefore gives
$s_n(A)=s_n(A^*)\le Ce^{-cn^{1/d}}$.
In particular, $A$ is compact and, after adjusting the constants,
$\lambda_n(A^*A)=s_n(A)^2\le Ce^{-cn^{1/d}}$.
If $\lambda_n(A^*A)\ge\tau$, then
$cn^{1/d}\le\log(C/\tau)$. Consequently,
\begin{equation}\label{eq:linear-count-upper}
N(\tau)
\le C'\bigl(\log(C/\tau)\bigr)^d
\le C''\bigl(\log(1/\tau)\bigr)^d
\end{equation}
for sufficiently small $\tau$.

\subsection{Lower bound for $N$}

Choose an axis-parallel cube $Q\Subset E$ of side length
$\ell>0$. For $n\in\mathbb Z^d$, put
$\kappa_n:=2\pi n/\ell$ and
$\varphi_n(x):=|Q|^{-1/2}e^{-i\kappa_n\cdot x}\mathbf1_Q(x)$.
The functions $(\varphi_n)_{n\in\mathbb Z^d}$ form an
orthonormal family in $L^2(E)$.

For each $\kappa\ne0$, choose $a\in\mathbb R^d$ with
$a\cdot\kappa=0$ and $|a|=|\kappa|/2$, which is possible
because $d\ge2$. Set
$\zeta_1:=a+i\kappa/2$ and $\zeta_2:=-a+i\kappa/2$.
Then
$\zeta_1\cdot\zeta_1
=|a|^2-|\kappa|^2/4+ia\cdot\kappa=0$,
and likewise $\zeta_2\cdot\zeta_2=0$.
Thus $U_j(x):=e^{\zeta_j\cdot x}$ are entire harmonic functions,
since $\Delta U_j=(\zeta_j\cdot\zeta_j)U_j$.
These are Calder\'on's exponential solutions
\cite{Calderon2006}. Moreover, $\zeta_1+\zeta_2=i\kappa$ and
$\zeta_1\cdot\zeta_2=-|a|^2-|\kappa|^2/4=-|\kappa|^2/2$,
so
\begin{equation}\label{eq:linear-harmonic-product}
\nabla U_1\cdot\nabla U_2
=
-\frac{|\kappa|^2}{2}e^{i\kappa\cdot x}.
\end{equation}

Let $F_j:=U_j|_{\partial\Omega}$. Boundedness of $\Omega$
and $|\zeta_j|\le C|\kappa|$ imply
$\|U_j\|_{H^1(\Omega)}
\le C(1+|\kappa|)e^{C|\kappa|}$.
The trace estimate therefore gives
\[
\|F_j\|_H
\le C\|U_j\|_{H^1(\Omega)}
\le C(1+|\kappa|)e^{C|\kappa|}
\le C'e^{C'|\kappa|}.
\]
Normalize the traces by
$f_\kappa:=F_1/\|F_1\|_H$ and
$g_\kappa:=F_2/\|F_2\|_H$.
Then $\|f_\kappa\|_H=\|g_\kappa\|_H=1$.
Since the harmonic extension of $F_j$ is exactly $U_j$,
\eqref{eq:linear-harmonic-product} yields
\[
\nabla u_{f_\kappa}\cdot\nabla u_{g_\kappa}
=
\alpha_\kappa e^{i\kappa\cdot x},
\qquad
\alpha_\kappa
:=
-\frac{|\kappa|^2}{2\|F_1\|_H\|F_2\|_H}.
\]
In particular,
$|\alpha_\kappa|\ge c|\kappa|^2e^{-C|\kappa|}$, and hence
$|\alpha_\kappa|\ge ce^{-CL}$ whenever $1\le|\kappa|\le L$.

For $L$ sufficiently large, define
$I_L:=\{n\in\mathbb Z^d:1\le|\kappa_n|\le L\}$ and
$M_L:=\#I_L$. Since $\kappa_n=2\pi n/\ell$, elementary
lattice counting gives $cL^d\le M_L\le CL^d$.
Let $V_L:=\operatorname{span}\{\varphi_n:n\in I_L\}$.
Then $\dim V_L=M_L$, and for
$q=\sum_{n\in I_L}c_n\varphi_n$ we have
$\|q\|_{L^2(E)}^2=\sum_{n\in I_L}|c_n|^2$.

For $m\in I_L$, the normalized exponential probes give
$(Aq)(f_{\kappa_m},g_{\kappa_m})
=\alpha_{\kappa_m}\int_Eq(x)e^{i\kappa_m\cdot x}\,dx$.
The integral extracts the corresponding Fourier coefficient:
\[
\int_E\varphi_n(x)e^{i\kappa_m\cdot x}\,dx
=
|Q|^{-1/2}\int_Qe^{i(\kappa_m-\kappa_n)\cdot x}\,dx
=
|Q|^{1/2}\delta_{nm}.
\]
Therefore
$(Aq)(f_{\kappa_m},g_{\kappa_m})
=\alpha_{\kappa_m}|Q|^{1/2}c_m$.
Squaring and summing over $m\in I_L$, we obtain
\begin{align}
\sum_{m\in I_L}
|(Aq)(f_{\kappa_m},g_{\kappa_m})|^2
&=
|Q|\sum_{m\in I_L}|\alpha_{\kappa_m}|^2|c_m|^2
\notag\\
&\ge
ce^{-CL}\sum_{m\in I_L}|c_m|^2
=
ce^{-CL}\|q\|_{L^2(E)}^2,
\label{eq:linear-probe-sum}
\end{align}
where the constants have been adjusted after squaring.

The probe pairs need not be orthogonal. What we use is that
each pair consists of unit vectors. For any $B\in\mathrm{HS}$,
Cauchy--Schwarz gives
$|B(f,g)|^2\le
\|B\|_{\mathrm{HS}}^2\|f\|_H^2\|g\|_H^2$.
Consequently,
\[
\sum_{m\in I_L}
|(Aq)(f_{\kappa_m},g_{\kappa_m})|^2
\le M_L\|Aq\|_{\mathrm{HS}}^2.
\]
Combining this with \eqref{eq:linear-probe-sum} and
$M_L\le CL^d$, we obtain
\[
\|Aq\|_{\mathrm{HS}}^2
\ge cM_L^{-1}e^{-CL}\|q\|_{L^2(E)}^2
\ge c'L^{-d}e^{-CL}\|q\|_{L^2(E)}^2.
\]
For $L\ge1$, the polynomial loss can be absorbed into the
exponential, since
$L^{-d}e^{-CL}=e^{-CL-d\log L}\ge e^{-(C+d)L}$.
Thus there are $c_0,C_0>0$ such that
\begin{equation}\label{eq:linear-subspace-lower}
\|Aq\|_{\mathrm{HS}}^2
\ge c_0e^{-C_0L}\|q\|_{L^2(E)}^2,
\qquad q\in V_L.
\end{equation}

Since $\langle A^*Aq,q\rangle=\|Aq\|_{\mathrm{HS}}^2$,
the min--max principle \cite{Kostenko} gives
\[
\lambda_{M_L}(A^*A)
\ge
\inf_{\substack{q\in V_L\\\|q\|_{L^2(E)}=1}}
\langle A^*Aq,q\rangle
\ge c_0e^{-C_0L}.
\]
Choose $L=(2C_0)^{-1}\log(1/\tau)$. For sufficiently small
$\tau$, this is large enough for the preceding estimates and
$c_0e^{-C_0L}=c_0\sqrt{\tau}\ge\tau$.
Hence
\begin{equation}\label{eq:linear-count-lower}
N(\tau)
\ge M_L
\ge cL^d
\ge c'\bigl(\log(1/\tau)\bigr)^d.
\end{equation}
Combining \eqref{eq:linear-count-upper} and
\eqref{eq:linear-count-lower} proves
Theorem~\ref{thm:linearized-counting}.

\section{Metric entropy of the nonlinear boundary data}

We now consider the full map $\gamma\mapsto\Lambda_\gamma$ with the
conductivity prescribed outside $E$. Fix $0<\lambda<L$ and a real
Lipschitz conductivity $\beta$ on $\Omega$ satisfying
$\lambda\leq\beta\leq L$. Let
\[
\mathcal K_\beta:=\bigl\{\gamma\in L^\infty(\Omega;\mathbb R):
  \lambda\leq\gamma\leq L\ \text{a.e.},\quad
  \gamma=\beta\ \text{a.e. on }\Omega\setminus E\bigr\}.
\]
Inside $E$ we assume only measurability and the stated ellipticity
bounds. Write $F_\beta(\gamma):=\Lambda_\gamma-\Lambda_\beta$.

For $\mathcal A\subset\mathrm{HS}$, let $N_\varepsilon(\mathcal A)$
be the smallest number of closed balls of radius $\varepsilon$ covering
$\mathcal A$. Its Kolmogorov $\varepsilon$-entropy is
$\mathcal H_\varepsilon(\mathcal A):=\log N_\varepsilon(\mathcal A)$
\cite{KT}. We estimate $\mathcal H_\varepsilon(F_\beta(\mathcal K_\beta))$
as $\varepsilon\downarrow0$.

\begin{theorem}\label{thm:entropy-upper}
The set $F_\beta(\mathcal K_\beta)$ is bounded in the Hilbert--Schmidt
norm. There exist $C>0$ and $\varepsilon_0\in(0,1)$, depending only
on $d,\Omega,E,\lambda,L$ and the chosen norm on $H$, such that
\[
\mathcal H_\varepsilon\bigl(F_\beta(\mathcal K_\beta)\bigr)
\le C\bigl(\log(1/\varepsilon)\bigr)^{2d-1},
\qquad 0<\varepsilon\le\varepsilon_0.
\]

The same bound holds if the covering centers are required to be
$F_\beta(\gamma_i)$ for conductivities $\gamma_i\in\mathcal K_\beta$.
\label{nonlinear upper bound}
\end{theorem}

\begin{remark}
It suffices that the prescribed exterior conductivity admit a
Lipschitz extension with values in $[\lambda,L]$. An exterior
function that is Lipschitz in the Euclidean distance has such an
extension by McShane's theorem~\cite{McShane}, followed by truncation to
$[\lambda,L]$. The extension inside $E$ does not change
$\mathcal K_\beta$; changing the reference DtN map translates the
data set and preserves its entropy.
\end{remark}

\subsection{The covering upper bound}

Set $T_\gamma f:=\nabla u_f^\gamma|_E$, viewed as an operator from
$H$ to $L^2(E;\mathbb C^d)$. The variational Dirichlet
estimate~\cite{RS} gives
\begin{equation}\label{eq:dirichlet-energy}
\|\nabla u_f^\gamma\|_{L^2(\Omega)}
\le C\|f\|_H,
\end{equation}
uniformly for $\gamma\in\mathcal K_\beta$.

\begin{lemma}\label{lem:singular-values}
The operator $T_\beta$ is compact, and there are $A,a>0$ such that
\[
s_j(T_\beta)\le A\exp\bigl(-a j^{1/(d-1)}\bigr),
\qquad j\ge1.
\]
The constants depend only on $d,\Omega,E,\lambda,L$ and the chosen
norm on $H$.
\end{lemma}

\begin{proof}
Put $\alpha:=1/(d-1)$. The local approximation theorem of Hauck, Ma
and Peterseim~\cite{HMP}, applied to $\beta I$ on
concentric cubes $Q$ and $3Q$, gives a space $V_{Q,n}$ of $\beta$-harmonic
functions on $Q$, of dimension at most $n$ modulo constants, such
that every $\beta$-harmonic function $v$ on $3Q$ satisfies
\begin{equation}\label{eq:local-approximation}
\inf_{z\in V_{Q,n}}
\left(\int_Q\beta|\nabla(v-z)|^2\,dx\right)^{1/2}
\le
C e^{-c n^\alpha}
\left(\int_{3Q}\beta|\nabla v|^2\,dx\right)^{1/2}.
\end{equation}
The constants depend only on dimension and ellipticity; increasing
$C$ covers the finitely many small values of $n$. The spaces and
projections are complexified for complex solutions.

Choose cubes $Q_1,\ldots,Q_J$ covering $\overline E$, with
$3Q_\ell\Subset\Omega$, and a measurable partition
$E=\bigsqcup_{\ell=1}^J E_\ell$ with $E_\ell\subset Q_\ell$.
Let $\Pi_{\ell,n}$ be the projection onto $V_{Q_\ell,n}$ for the
energy inner product
$\int_{Q_\ell}\beta\nabla v\cdot\overline{\nabla z}\,dx$,
with constants identified, and define
\[
(R_nf)|_{E_\ell}
:=
\nabla\Pi_{\ell,n}(u_f^\beta|_{Q_\ell})|_{E_\ell}.
\]
Then $\operatorname{rank}R_n\le Jn$. By
\eqref{eq:local-approximation} and \eqref{eq:dirichlet-energy},
\begin{align*}
\|(T_\beta-R_n)f\|_{L^2(E)}^2
&\le \lambda^{-1}\sum_{\ell=1}^J
\int_{Q_\ell}\beta
\bigl|\nabla u_f^\beta-
\nabla\Pi_{\ell,n}(u_f^\beta|_{Q_\ell})\bigr|^2\,dx\\
&\le C e^{-2c n^\alpha}
\sum_{\ell=1}^J\int_{3Q_\ell}\beta|\nabla u_f^\beta|^2\,dx\\
&\le C' e^{-2c n^\alpha}\|f\|_H^2.
\end{align*}
Thus $T_\beta$ is compact, and the characterization by finite-rank
approximation~\cite{Kostenko} gives
\[
s_{Jn+1}(T_\beta)\le\|T_\beta-R_n\|
\le C e^{-c n^\alpha}.
\]
Monotonicity of the singular values, with the fixed number $J$
absorbed into the constants, proves the claim.
\end{proof}

\begin{lemma}\label{lem:nonlinear-domination}
For every $\gamma\in\mathcal K_\beta$ and $f\in H$,
\begin{equation}\label{eq:gradient-domination}
\|T_\gamma f\|_{L^2(E)}
\le\frac L\lambda\|T_\beta f\|_{L^2(E)}.
\end{equation}
Consequently,
\begin{equation}\label{eq:form-domination}
|F_\beta(\gamma)(f,g)|
\le M\|T_\beta f\|_{L^2(E)}\|T_\beta g\|_{L^2(E)},
\qquad M:=(L-\lambda)\frac L\lambda.
\end{equation}
\end{lemma}

\begin{proof}
Let $w:=u_f^\gamma-u_f^\beta$. Then $w\in H_0^1(\Omega)$ and
\[
\int_\Omega\gamma\nabla w\cdot\nabla\phi\,dx
=-\int_E(\gamma-\beta)\nabla u_f^\beta\cdot\nabla\phi\,dx,
\qquad \phi\in H_0^1(\Omega).
\]
Taking $\phi=\overline w$, we obtain
\[
\lambda\|\nabla w\|_{L^2(\Omega)}^2
\le\|\gamma-\beta\|_\infty
\|T_\beta f\|_{L^2(E)}\|\nabla w\|_{L^2(\Omega)}.
\]
Hence
\[
\|T_\gamma f\|_{L^2(E)}
\le
\left(1+\frac{\|\gamma-\beta\|_\infty}{\lambda}\right)
\|T_\beta f\|_{L^2(E)}
\le\frac L\lambda\|T_\beta f\|_{L^2(E)}.
\]
Since $u_f^\gamma-u_f^\beta$ and $u_g^\gamma-u_g^\beta$ belong to
$H_0^1(\Omega)$, the weak equations give Alessandrini's
identity~\cite{Alessandrini1988} (see also
\cite[Eq.~(2.4)]{ADHG2017}),
\[
F_\beta(\gamma)(f,g)
=\int_E(\gamma-\beta)T_\gamma f\cdot T_\beta g\,dx.
\]
The second estimate now follows from
\eqref{eq:gradient-domination} and Cauchy--Schwarz.
\end{proof}

\begin{proof}[Proof of Theorem~\ref{thm:entropy-upper}]
Write $s_j:=s_j(T_\beta)$ and $\alpha:=1/(d-1)$. Since $T_\beta$
commutes with conjugation, choose real-valued right singular vectors
$(e_j)$ forming an orthonormal basis of $(\ker T_\beta)^\perp$. Estimate
\eqref{eq:form-domination} makes every $F_\beta(\gamma)$ vanish
when either argument lies in $\ker T_\beta$, so the kernel
contributes no matrix entries. In the finite-rank case, matrices
and singular-value sequences are padded with zeros, and projections
use only the available singular vectors.

By \eqref{eq:form-domination},
$|F_\beta(\gamma)(e_j,e_k)|\le Ms_js_k$. The sum
$S:=\sum_{j\ge1}s_j^2$ is uniformly bounded by
Lemma~\ref{lem:singular-values}; hence
\begin{equation}\label{eq:hs-radius}
\|F_\beta(\gamma)\|_{\mathrm{HS}}^2
\le M^2\sum_{j,k\ge1}s_j^2s_k^2=M^2S^2.
\end{equation}

Let $P_m$ be the orthogonal projection onto the first $m$ right
singular vectors and set
\[
F_\beta^{(m)}(\gamma)(f,g)
:=F_\beta(\gamma)(P_mf,P_mg).
\]
Writing $S_m:=\sum_{j>m}s_j^2$, we obtain
\begin{align*}
\|F_\beta(\gamma)-F_\beta^{(m)}(\gamma)\|_{\mathrm{HS}}^2
&=\sum_{\substack{j,k\ge1\\j>m\ \mathrm{or}\ k>m}}
|F_\beta(\gamma)(e_j,e_k)|^2\\
&\le M^2\bigl(S^2-(S-S_m)^2\bigr)\\
&\le2M^2SS_m.
\end{align*}
Moreover,
\[
S_m\le A^2\sum_{j>m}e^{-2aj^\alpha}
\le A^2e^{-am^\alpha}\sum_{j\ge1}e^{-aj^\alpha}.
\]
Consequently, there are uniform constants $B\ge1$ and $b>0$ such that
\begin{equation}\label{eq:uniform-truncation}
\sup_{\gamma\in\mathcal K_\beta}
\|F_\beta(\gamma)-F_\beta^{(m)}(\gamma)\|_{\mathrm{HS}}
\le B e^{-bm^{1/(d-1)}}.
\end{equation}

Fix $0<\varepsilon<1$ and choose
\begin{equation}\label{eq:choice-m}
m:=\left\lceil
\left(\frac1b\log\frac{4B}{\varepsilon}\right)^{d-1}
\right\rceil.
\end{equation}
Then the truncation error is at most $\varepsilon/4$. The truncated
forms are real symmetric $m\times m$ matrices, with the Frobenius
norm, and lie in a ball of radius $R:=MS$ by \eqref{eq:hs-radius}.
This matrix space has real dimension $D_m:=m(m+1)/2\le m^2$.

Choose a maximal $\varepsilon/2$-separated family
$\{F_\beta^{(m)}(\gamma_i)\}_{i=1}^N$ in the truncated image.
The balls of radius $\varepsilon/4$ around these points are
disjoint, up to their boundaries, and are contained in the ball of
radius $R+\varepsilon/4$. The usual volume comparison
\cite[Proposition~4.2.12]{Vershynin} gives
\begin{equation}\label{eq:finite-cover-count}
N\le\left(1+\frac{4R}{\varepsilon}\right)^{D_m}.
\end{equation}
In particular, a maximal separated family can be chosen finite.
By maximality, each $\gamma\in\mathcal K_\beta$ has a representative
$\gamma_i$ with
\[
\|F_\beta^{(m)}(\gamma)-F_\beta^{(m)}(\gamma_i)\|_{\mathrm{HS}}
\le\varepsilon/2.
\]
Adding the two truncation errors gives
\[
\|F_\beta(\gamma)-F_\beta(\gamma_i)\|_{\mathrm{HS}}
\le\varepsilon/4+\varepsilon/2+\varepsilon/4=\varepsilon.
\]
Thus the covering centers come from admissible conductivities.

Finally, \eqref{eq:choice-m} and \eqref{eq:finite-cover-count} imply
\[
\mathcal H_\varepsilon\bigl(F_\beta(\mathcal K_\beta)\bigr)
\le m^2\log\left(1+\frac{4R}{\varepsilon}\right)
\le C_0\left(\log\frac{C_1}{\varepsilon}\right)^{2d-1}.
\]
For sufficiently small $\varepsilon_0$, this is bounded by
$C(\log(1/\varepsilon))^{2d-1}$ for
$0<\varepsilon\le\varepsilon_0$.
\end{proof}

\begin{corollary}\label{entropy upper bound corollary}
Let
\[
\mathcal D_\beta
:=\{\Lambda_{\gamma_1}-\Lambda_{\gamma_2}:
\gamma_1,\gamma_2\in\mathcal K_\beta\}.
\]
Then, after adjusting the constants,
\[
\mathcal H_\varepsilon(\mathcal D_\beta)
\le C\bigl(\log(1/\varepsilon)\bigr)^{2d-1},
\qquad 0<\varepsilon\le\varepsilon_0.
\]
\end{corollary}

\begin{proof}
Taking differences of the centers of an $\varepsilon/2$-cover of
$F_\beta(\mathcal K_\beta)$ gives an $\varepsilon$-cover of
$\mathcal D_\beta$. Hence
\[
N_\varepsilon(\mathcal D_\beta)
\le N_{\varepsilon/2}\bigl(F_\beta(\mathcal K_\beta)\bigr)^2,
\]
and the claim follows from Theorem~\ref{thm:entropy-upper}.
\end{proof}

\begin{remark}
The upper bound uses only boundedness and uniform ellipticity of
$\beta$. In particular, its constants are independent of the
Lipschitz seminorm. We retain the Lipschitz assumption to use the
same class of exterior conductivities in the lower-bound argument.
\end{remark}

\subsection{Localized Fourier probes}

\begin{theorem}\label{thm:entropy-lower}
There exist $c>0$ and $\varepsilon_0\in(0,1)$ such that
\[
\mathcal H_\varepsilon\bigl(F_\beta(\mathcal K_\beta)\bigr)
\ge c\bigl(\log(1/\varepsilon)\bigr)^{d+1},
\qquad 0<\varepsilon\le\varepsilon_0.
\]
The constants may depend on $d,\Omega,E,\lambda,L$, the chosen
norm on $H$, and the Lipschitz norm of $\beta$.
\label{nonlinear lower bound}
\end{theorem}

Since the lower bound only requires a family of examples in $\mathcal K_\beta$, we may choose a reference conductivity that is constant on a ball inside $E$. On that ball the harmonic exponentials from Section 2 are again solutions of the conductivity equation. Choose a cube $Q$ and concentric balls $B_1,B_2$ with
$Q\Subset B_1\Subset B_2\Subset E$, and take
$\chi\in C_c^\infty(E)$ with $0\le\chi\le1$ and $\chi=1$
on a neighborhood of $\overline{B_2}$. Set
$a:=(\lambda+L)/2$ and
$\gamma_*:=(1-\chi)\beta+\chi a$.
Then $\gamma_*\in\mathcal K_\beta$ is Lipschitz and equals $a$
on $B_2$. We consider real perturbations $q$ supported in $Q$,
extended by zero to $\Omega$. In particular,
$\gamma_*+q\in\mathcal K_\beta$ whenever
$\|q\|_{L^\infty(Q)}\le(L-\lambda)/2$.

Write $T_*:=T_{\gamma_*}$ and define the linearized form by
\[
A_*q(f,g)
:=\int_Q q\,\nabla u_f^{\gamma_*}\cdot
                 \nabla u_g^{\gamma_*}\,dx.
\]
By Lemma~\ref{lem:singular-values}, applied with background
$\gamma_*$, the operator $T_*$ is Hilbert--Schmidt. Consequently,
$A_*q$ is a Hilbert--Schmidt form for every $q\in L^\infty(Q)$.

We next record the finite-dimensional estimates needed from the
linear argument. The $L^\infty$ estimate ensures that the
perturbations used below remain admissible.

\begin{lemma}\label{lem:nonlinear-linear-input}
There exist constants $c_0,c_1,C_0,C_1>0$ and $t_0\ge1$ such
that, for every $t\ge t_0$, there is a real subspace
$V_t\subset L^\infty(Q;\mathbb R)$ satisfying
\[
\begin{gathered}
M_t:=\dim_{\mathbb R}V_t,
\qquad c_1t^d\le M_t\le C_1t^d,\\
\|q\|_{L^\infty(Q)}
\le C_1t^{d/2}\|q\|_{L^2(Q)},
\qquad
\|A_*q\|_{\mathrm{HS}}
\ge c_0e^{-C_0t}\|q\|_{L^2(Q)}
\quad(q\in V_t).
\end{gathered}
\]
\end{lemma}

\begin{proof}
Let $\ell$ be the side length of $Q$, put
$\kappa_n:=2\pi n/\ell$, and take
$I_t:=\{n\in\mathbb Z^d:1\le|\kappa_n|\le t\}$.
Using the same Fourier subspaces as in the linear proof, set
\[
V_t:=
\left\{
q=\sum_{n\in I_t}c_n|Q|^{-1/2}e^{-i\kappa_n\cdot x}:
c_{-n}=\overline{c_n}
\right\}.
\]
The symmetry condition makes $q$ real-valued. Since $I_t$ is
symmetric, $\dim_{\mathbb R}V_t=\#I_t\asymp t^d$.
Orthogonality and Cauchy--Schwarz give
\[
\|q\|_{L^2(Q)}^2=\sum_{n\in I_t}|c_n|^2,
\qquad
\|q\|_{L^\infty(Q)}
\le |Q|^{-1/2}\sum_{n\in I_t}|c_n|
\le C\sqrt{M_t}\,\|q\|_{L^2(Q)}.
\]

For each $n\in I_t$, choose a real vector $b_n$ with
$b_n\perp\kappa_n$ and $|b_n|=|\kappa_n|/2$, and set
$\zeta_{n,1}:=b_n+i\kappa_n/2$,
$\zeta_{n,2}:=-b_n+i\kappa_n/2$.
The functions $h_{n,j}(x):=e^{\zeta_{n,j}\cdot x}$ are harmonic
and satisfy
\[
\nabla h_{n,1}\cdot\nabla h_{n,2}
=-\frac{|\kappa_n|^2}{2}e^{i\kappa_n\cdot x},
\qquad
\|h_{n,j}\|_{H^1(B_2)}\le Ce^{Ct}.
\]
Since $\gamma_*=a$ on $B_2$, these functions also solve the
$\gamma_*$-conductivity equation there.

To obtain global solutions, apply quantitative Runge
approximation~\cite{RS} with inner domain $B_1$,
intermediate domain $B_2$, and outer domain $\Omega$.
Here $\Omega\setminus\overline{B_1}$ is connected, and the
accessible boundary is all of $\partial\Omega$.
For every $\eta\in(0,1)$ there are global
$\gamma_*$-harmonic functions $U_{n,j}$, with boundary traces
$f_{n,j}$, such that
\[
\|U_{n,j}-h_{n,j}\|_{L^2(B_1)}
\le \eta\|h_{n,j}\|_{H^1(B_2)},
\qquad
\|f_{n,j}\|_H
\le C\eta^{-\mu}\|h_{n,j}\|_{L^2(B_1)}
\le C\eta^{-\mu}e^{Ct}.
\]
The result applies to complex functions by treating real and
imaginary parts separately. Moreover, $U_{n,j}-h_{n,j}$ is
harmonic on $B_1$, so the interior Caccioppoli estimate gives
$$
\|\nabla(U_{n,j}-h_{n,j})\|_{L^2(Q)}
\le C\|U_{n,j}-h_{n,j}\|_{L^2(B_1)}
\le C\eta e^{Ct}.
$$

For $q\in V_t$, Fourier orthogonality gives
$\int_Qq\,e^{i\kappa_n\cdot x}\,dx=|Q|^{1/2}c_n$.
Expanding the difference of the two gradient products and
using Cauchy--Schwarz therefore yields
\[
\left|
A_*q(f_{n,1},f_{n,2})
+\frac{|\kappa_n|^2}{2}|Q|^{1/2}c_n
\right|
\le C\sqrt{M_t}\,\eta e^{Ct}\|q\|_{L^2(Q)}.
\]
Taking the Euclidean norm over $n\in I_t$, and using
$|\kappa_n|\ge1$, we obtain
\[
\left(\sum_{n\in I_t}
|A_*q(f_{n,1},f_{n,2})|^2\right)^{1/2}
\ge
\left(\frac{|Q|^{1/2}}2-CM_t\eta e^{Ct}\right)
\|q\|_{L^2(Q)}.
\]
Choose $\eta=\eta_t:=c_*M_t^{-1}e^{-C_*t}$, with $c_*>0$
sufficiently small and $C_*>0$ sufficiently large.
The coefficient on the right is then bounded below by a
positive constant. The Runge estimate still gives
$\|f_{n,j}\|_H\le Ce^{Ct}$, since $M_t$ grows polynomially.

Finally, the Hilbert--Schmidt evaluation bound gives
\[
c\|q\|_{L^2(Q)}
\le
\left(\sum_{n\in I_t}
|A_*q(f_{n,1},f_{n,2})|^2\right)^{1/2}
\le C\sqrt{M_t}\,e^{Ct}\|A_*q\|_{\mathrm{HS}}.
\]
Absorbing the polynomial factor into the exponential proves
the claimed lower bound.
\end{proof}

\subsection{The packing lower bound}

\begin{lemma}\label{lem:quadratic-remainder}
For every admissible perturbation $q$,
\[
F_\beta(\gamma_*+q)-F_\beta(\gamma_*)
=A_*q+\mathcal R(q),
\qquad
\|\mathcal R(q)\|_{\mathrm{HS}}
\le C_R\|q\|_{L^\infty(Q)}^2,
\]
where $C_R$ is independent of $q$.
\end{lemma}

\begin{proof}
Let $w_f:=u_f^{\gamma_*+q}-u_f^{\gamma_*}$.
Since the boundary data agree, $w_f\in H_0^1(\Omega)$, and
\[
\int_\Omega(\gamma_*+q)\nabla w_f\cdot\nabla\phi\,dx
=-\int_Q q\,\nabla u_f^{\gamma_*}\cdot\nabla\phi\,dx,
\qquad \phi\in H_0^1(\Omega).
\]
Taking $\phi=\overline{w_f}$ and using ellipticity gives
\[
\lambda\|\nabla w_f\|_{L^2(\Omega)}^2
\le
\|q\|_{L^\infty(Q)}
\|T_*f\|_{L^2(E)}
\|\nabla w_f\|_{L^2(\Omega)}.
\]
Hence
$\|\nabla w_f\|_{L^2(\Omega)}
\le\lambda^{-1}\|q\|_{L^\infty(Q)}\|T_*f\|_{L^2(E)}$.

Alessandrini's identity~\cite{Alessandrini1988} yields
\begin{align*}
\bigl(F_\beta(\gamma_*+q)-F_\beta(\gamma_*)\bigr)(f,g)
&=\int_Q q\,\nabla u_f^{\gamma_*+q}
                  \cdot\nabla u_g^{\gamma_*}\,dx\\
&=A_*q(f,g)+
  \int_Q q\,\nabla w_f\cdot\nabla u_g^{\gamma_*}\,dx.
\end{align*}
Thus
\[
|\mathcal R(q)(f,g)|
\le\lambda^{-1}\|q\|_{L^\infty(Q)}^2
   \|T_*f\|_{L^2(E)}\|T_*g\|_{L^2(E)}.
\]
For an orthonormal basis $(e_j)$ of $H$, summing the squared
matrix entries gives
\[
\|\mathcal R(q)\|_{\mathrm{HS}}^2
\le
\lambda^{-2}\|q\|_{L^\infty(Q)}^4
\left(\sum_j\|T_*e_j\|_{L^2(E)}^2\right)^2.
\]
The conclusion follows with
$C_R:=\lambda^{-1}\|T_*\|_{\mathrm{HS}}^2$.
\end{proof}

\begin{proof}[Proof of Theorem~\ref{thm:entropy-lower}]
Choose $r_0>0$ such that
$r_0\le(L-\lambda)/2$ and $2C_Rr_0^2\le1$.
For sufficiently small $\varepsilon>0$, set
\[
\begin{aligned}
t&:=\frac{1}{4C_0}\log(1/\varepsilon),
&\qquad \sigma_t&:=c_0e^{-C_0t}=c_0\varepsilon^{1/4},\\
\rho&:=\frac{r_0\sqrt{\varepsilon}}{C_1t^{d/2}},
&\qquad \delta&:=\frac{4\varepsilon}{\sigma_t}.
\end{aligned}
\]
We take $\varepsilon$ small enough that $t\ge t_0$.

Consider the $L^2(Q)$-ball
$\mathcal B_t(\rho):=\{q\in V_t:\|q\|_{L^2(Q)}\le\rho\}$.
By Lemma~\ref{lem:nonlinear-linear-input}, every
$q\in\mathcal B_t(\rho)$ satisfies
$\|q\|_{L^\infty(Q)}\le r_0\sqrt{\varepsilon}$.
Thus $\gamma_*+q\in\mathcal K_\beta$, and
Lemma~\ref{lem:quadratic-remainder} gives
\[
\sup_{q\in\mathcal B_t(\rho)}
\|\mathcal R(q)\|_{\mathrm{HS}}
\le C_Rr_0^2\varepsilon.
\]

Choose a maximal finite family $q_1,\ldots,q_M$ in
$\mathcal B_t(\rho)$ with
$\|q_i-q_j\|_{L^2(Q)}>\delta$ for $i\ne j$.
By maximality, the balls of radius $\delta$ centered at these
points cover $\mathcal B_t(\rho)$.
Volume comparison in the real $M_t$-dimensional space $V_t$
therefore gives
\[
M\ge\left(\frac{\rho}{\delta}\right)^{M_t}.
\]
Set $\gamma_i:=\gamma_*+q_i$. For distinct $i,j$, the linear
lower bound and the quadratic remainder estimate imply
\begin{align*}
\|F_\beta(\gamma_i)-F_\beta(\gamma_j)\|_{\mathrm{HS}}
&\ge
\|A_*(q_i-q_j)\|_{\mathrm{HS}}
-\|\mathcal R(q_i)\|_{\mathrm{HS}}
-\|\mathcal R(q_j)\|_{\mathrm{HS}}\\
&\ge
\sigma_t\|q_i-q_j\|_{L^2(Q)}
-2C_Rr_0^2\varepsilon\\
&>4\varepsilon-2C_Rr_0^2\varepsilon
\ge3\varepsilon.
\end{align*}
A Hilbert--Schmidt ball of radius $\varepsilon$ can therefore
contain at most one of these images. Consequently,
\[
N_\varepsilon\bigl(F_\beta(\mathcal K_\beta)\bigr)
\ge M
\ge\left(\frac{\rho}{\delta}\right)^{M_t}.
\]

By our choices of $\rho,\delta$ and $\sigma_t$,
\[
\frac{\rho}{\delta}
=\frac{r_0c_0}{4C_1}\,
  \varepsilon^{-1/4}t^{-d/2},
\qquad
\log\frac{\rho}{\delta}
=\frac14\log(1/\varepsilon)
-\frac d2\log t
+\log\frac{r_0c_0}{4C_1}.
\]
Since $t=(4C_0)^{-1}\log(1/\varepsilon)$, for sufficiently
small $\varepsilon$ the last expression is at least
$\frac18\log(1/\varepsilon)$. Hence
\[
\mathcal H_\varepsilon
\bigl(F_\beta(\mathcal K_\beta)\bigr)
\ge M_t\log\frac{\rho}{\delta}
\ge \frac{c_1}{8}t^d\log(1/\varepsilon)
\ge c\bigl(\log(1/\varepsilon)\bigr)^{d+1}.
\]
Although the perturbation radius depends on $\varepsilon$,
all the conductivities used in these packings belong to the
same fixed class $\mathcal K_\beta$.
\end{proof}

Since $\beta\in\mathcal K_\beta$, the set $F_\beta(\mathcal K_\beta)$
is contained in $\mathcal D_\beta$. The lower entropy bound therefore
holds for $\mathcal D_\beta$ as well. Together with the upper bound
for this set, this gives
$\mathcal H_\varepsilon(\mathcal D_\beta)\asymp
(\log(1/\varepsilon))^3$ in dimension two.

\begin{remark}[Nonlinear entropy in the operator norm]
For a bounded bilinear form $B$ on $H\times H$, write
$\|B\|_{\mathrm{op}}=\sup_{\|f\|_H=\|g\|_H=1}|B(f,g)|$.
For a difference of Dirichlet-to-Neumann maps this is the operator
norm from $H^{1/2}(\partial\Omega)$ to $H^{-1/2}(\partial\Omega)$,
with the induced dual norm. Denote the entropy in this norm by
$\mathcal H^{\mathrm{op}}_\varepsilon$. Since
$\|B\|_{\mathrm{op}}\leq\|B\|_{\mathrm{HS}}$, the covers in Theorem~\ref{thm:entropy-upper}
and Corollary~\ref{entropy upper bound corollary} give the same upper bounds in operator norm.

For the lower bound, the boundary data in the proof of Lemma~\ref{lem:nonlinear-linear-input} give,
for $q\in V_t$,
\[
c\|q\|_{L^2(Q)}
 \leq \left(\sum_{n\in I_t}|A_*q(f_{n,1},f_{n,2})|^2\right)^{1/2}
 \leq C\sqrt{M_t}\,e^{Ct}\|A_*q\|_{\mathrm{op}} .
\]

The last inequality follows by bounding each evaluation by
$\|A_*q\|_{\mathrm{op}}\|f_{n,1}\|_H\|f_{n,2}\|_H$.
As $M_t\asymp t^d$, we obtain
$\|A_*q\|_{\mathrm{op}}\geq c_0e^{-C_0t}\|q\|_{L^2(Q)}$ after
adjusting the constants. Lemma~\ref{lem:quadratic-remainder} also gives
$\|\mathcal R(q)\|_{\mathrm{op}}\leq\|\mathcal R(q)\|_{\mathrm{HS}}
\leq C_R\|q\|_{L^\infty(Q)}^2$. With these constants the packing in
the proof of Theorem~\ref{thm:entropy-lower} has linearized responses more than
$4\varepsilon$ apart in operator norm, and each remainder has norm at
most $\varepsilon/2$. The full responses are therefore more than
$3\varepsilon$ apart. It follows that, for all sufficiently small
$\varepsilon>0$,
\[
c\bigl(\log(1/\varepsilon)\bigr)^{d+1}
 \leq \mathcal H^{\mathrm{op}}_\varepsilon
       \bigl(F_\beta(\mathcal K_\beta)\bigr)
 \leq C\bigl(\log(1/\varepsilon)\bigr)^{2d-1}.
\]
The same bounds hold for $\mathcal D_\beta$, since it contains the
same packing. In dimension two the entropy of either set thus has
order $\bigl(\log(1/\varepsilon)\bigr)^3$ in operator norm.
\end{remark}

\section{Discussion}

At a boundary noise level $\varepsilon$, Theorem~\ref{thm:linearized-counting} gives the number of
singular values of the linearization $A$ at conductivity one that are
at least $\varepsilon$. In the notation of Section~2 this number is
$N(\varepsilon^2)$, since the eigenvalues of $A^*A$ are $s_j(A)^2$.
Writing $\ell_\varepsilon=\log(1/\varepsilon)$, we have
\[
N(\varepsilon^2)=\#\{j:s_j(A)\geq\varepsilon\}
   \asymp \ell_\varepsilon^d.
\]
This counts each singular direction once. To compare with the nonlinear
result, let $B$ be the unit ball of $L^2(E)$ and consider all linearized
responses $Aq$ with $q\in B$. Here $N_\varepsilon(A(B))$ counts the
Hilbert--Schmidt balls of radius $\varepsilon$ needed to cover these
responses. Its logarithm satisfies
\[
\mathcal H_\varepsilon(A(B))=\log N_\varepsilon(A(B))
   \asymp \ell_\varepsilon^{d+1}.
\]

The extra power of $\ell_\varepsilon$ can be seen from the ellipsoid
$A(B)$, whose semiaxes have lengths $s_j(A)$. By Theorem~\ref{thm:linearized-counting}, the number
$m$ of semiaxes of length at least $\sqrt{\varepsilon}$ is of order
$\ell_\varepsilon^d$. In their span the ellipsoid contains a ball of
radius $\sqrt{\varepsilon}$. By volume comparison, covering this ball
at radius $\varepsilon$ takes at least $(1/\sqrt{\varepsilon})^m$
balls. Taking logarithms gives the lower bound
$\tfrac12 m\ell_\varepsilon$. For the upper
bound, retaining the semiaxes of length at least $\varepsilon/2$ gives
a truncation error of at most $\varepsilon/2$ on $B$. The truncated
ellipsoid has dimension $O(\ell_\varepsilon^d)$ and radius at most
$\|A\|$. A volume estimate at radius $\varepsilon/2$ therefore gives
$\log N_\varepsilon(A(B))\leq C\ell_\varepsilon^{d+1}$.
See also~\cite[Lemma~3.9]{KochRulandSalo2025}.

For the full conductivity map, Theorems~\ref{thm:entropy-upper} and ~\ref{thm:entropy-lower} give
$$
c\,\ell_\varepsilon^{d+1}
 \leq \log N_\varepsilon\bigl(F_\beta(\mathcal K_\beta)\bigr)
 \leq C\,\ell_\varepsilon^{2d-1}.
$$
For the lower bound we perturb an admissible reference conductivity.
The linearized responses of the chosen perturbations are more than
$4\varepsilon$ apart, and the quadratic remainder for each has norm at most
$\varepsilon/2$. Their full boundary responses are therefore more than
$3\varepsilon$ apart.
In dimension two, both exponents are $3$. Thus the linearized and
nonlinear images both have entropy of order $\ell_\varepsilon^3$.
For either image the number of covering balls lies between
$\exp(c\ell_\varepsilon^3)$ and $\exp(C\ell_\varepsilon^3)$, whereas
the number of linear singular directions above the noise level is of
order $\ell_\varepsilon^2$. In dimensions $d\geq3$, our bounds leave
open whether the nonlinear image has the same entropy order as the
linearized image.

\section*{Acknowledgements}
We thank Guy David, Max Engelstein, and Alberto Pacati for useful discussions. RJ and LS  acknowledge support from the Simons Foundation (00017375, RJ).
The work of RJ and LS was supported in part by project PID2022-141125NB-I00
and by the ``Center of Excellence Mar\'ia de Maeztu 2025--2029'' award
to the ICCUB, funded by grant CEX2024-001451-M from
AEI/10.13039/501100011033.

\section*{AI Usage}

We used AI tools in deriving some of the results in Section~3,
particularly in the proof of Lemma~\ref{lem:nonlinear-linear-input}. We also used these tools
to improve the grammar in parts of the manuscript. 

\bibliographystyle{plain}
\bibliography{Bibliography}

@article{Alessandrini1988,
  author  = {Alessandrini, Giovanni},
  title   = {Stable determination of conductivity by boundary measurements},
  journal = {Applicable Analysis},
  volume  = {27},
  number  = {1--3},
  pages   = {153--172},
  year    = {1988},
  doi     = {10.1080/00036818808839730}
}

@article{ADHG2017,
  author  = {Alessandrini, Giovanni and de Hoop, Maarten V. and
             Gaburro, Romina},
  title   = {Uniqueness for the electrostatic inverse boundary value
             problem with piecewise constant anisotropic conductivities},
  journal = {Inverse Problems},
  volume  = {33},
  number  = {12},
  pages   = {125013},
  year    = {2017},
  doi     = {10.1088/1361-6420/aa982d}
}

@misc{HMP,
  author        = {Hauck, Moritz and Ma, Chupeng and Peterseim, Daniel},
  title         = {Local Approximation of {$A$}-Harmonic Functions
                   and the Super-Localization Conjecture},
  year          = {2026},
  eprint        = {2610.25230},
  archivePrefix = {arXiv},
  primaryClass  = {math.NA},
  doi           = {10.48550/arXiv.2610.25230},
  url           = {https://arxiv.org/abs/2610.25230}
}

@article{KT,
  author  = {Kolmogorov, A. N. and Tikhomirov, V. M.},
  title   = {{$\varepsilon$}-entropy and {$\varepsilon$}-capacity
             of sets in function spaces},
  journal = {Uspekhi Matematicheskikh Nauk},
  volume  = {14},
  number  = {2(86)},
  pages   = {3--86},
  year    = {1959},
  url     = {https://www.mathnet.ru/eng/rm7289}
}

@misc{Kostenko,
  author = {Kostenko, Aleksey},
  title  = {Trace Ideals with Applications},
  note   = {Lecture notes},
  url    = {https://users.fmf.uni-lj.si/kostenko/teach/IdealsNotes.pdf}
}

@article{McShane,
  author  = {McShane, E. J.},
  title   = {Extension of range of functions},
  journal = {Bulletin of the American Mathematical Society},
  volume  = {40},
  number  = {12},
  pages   = {837--842},
  year    = {1934},
  doi     = {10.1090/S0002-9904-1934-05978-0}
}

@article{RS,
  author  = {R{\"u}land, Angkana and Salo, Mikko},
  title   = {Quantitative {Runge} Approximation and Inverse Problems},
  journal = {International Mathematics Research Notices},
  volume  = {2019},
  number  = {20},
  pages   = {6216--6234},
  year    = {2019},
  doi     = {10.1093/imrn/rnx301}
}

@book{Vershynin,
  author    = {Vershynin, Roman},
  title     = {High-Dimensional Probability:
               An Introduction with Applications in Data Science},
  series    = {Cambridge Series in Statistical and Probabilistic Mathematics},
  volume    = {47},
  publisher = {Cambridge University Press},
  year      = {2018},
  doi       = {10.1017/9781108231596}
}

@book{Evans2010,
  author    = {Evans, Lawrence C.},
  title     = {Partial Differential Equations},
  edition   = {2},
  series    = {Graduate Studies in Mathematics},
  volume    = {19},
  publisher = {American Mathematical Society},
  address   = {Providence, RI},
  year      = {2010},
  isbn      = {978-0-8218-4974-3}
}

@article{Calderon2006,
  author  = {Calder{\'o}n, Alberto P.},
  title   = {On an inverse boundary value problem},
  journal = {Computational and Applied Mathematics},
  volume  = {25},
  number  = {2--3},
  pages   = {133--138},
  year    = {2006},
  doi     = {10.1590/S0101-82052006000200002},
  note    = {Reprint of the 1980 original}
}

@article{AbrahamNickl2019,
  author  = {Abraham, Kweku and Nickl, Richard},
  title   = {On statistical {C}alder{\'o}n problems},
  journal = {Mathematical Statistics and Learning},
  volume  = {2},
  number  = {2},
  pages   = {165--216},
  year    = {2019}
}

@article{AdlerYoumaranLionheart2008,
  author  = {Adler, Andy and Youmaran, Richard and Lionheart, William R. B.},
  title   = {A measure of the information content of {EIT} data},
  journal = {Physiological Measurement},
  volume  = {29},
  number  = {6},
  pages   = {S101--S109},
  year    = {2008}
}

@article{Alessandrini2007,
  author  = {Alessandrini, Giovanni},
  title   = {Open issues of stability for the inverse conductivity problem},
  journal = {Journal of Inverse and Ill-Posed Problems},
  volume  = {15},
  number  = {5},
  pages   = {451--460},
  year    = {2007}
}

@article{AlessandriniScapin2017,
  author  = {Alessandrini, Giovanni and Scapin, Andrea},
  title   = {Depth dependent resolution in electrical impedance tomography},
  journal = {Journal of Inverse and Ill-Posed Problems},
  volume  = {25},
  number  = {3},
  pages   = {391--402},
  year    = {2017}
}

@article{AlessandriniVessella2005,
  author  = {Alessandrini, Giovanni and Vessella, Sergio},
  title   = {Lipschitz stability for the inverse conductivity problem},
  journal = {Advances in Applied Mathematics},
  volume  = {35},
  number  = {2},
  pages   = {207--241},
  year    = {2005}
}

@article{AllersSantosa1991,
  author  = {Allers, A. and Santosa, Fadil},
  title   = {Stability and resolution analysis of a linearized problem in electrical impedance tomography},
  journal = {Inverse Problems},
  volume  = {7},
  number  = {4},
  pages   = {515--533},
  year    = {1991},
  doi     = {10.1088/0266-5611/7/4/003}
}

@article{AstalaPaivarinta2006,
  author  = {Astala, Kari and P{\"a}iv{\"a}rinta, Lassi},
  title   = {Calder{\'o}n's inverse conductivity problem in the plane},
  journal = {Annals of Mathematics (2)},
  volume  = {163},
  number  = {1},
  pages   = {265--299},
  year    = {2006}
}

@article{BischHirvensaloHyvonen2025,
  author  = {Bisch, Joanna and Hirvensalo, Markus and Hyv{\"o}nen, Nuutti},
  title   = {Continuity of the linearized forward map of electrical impedance tomography from square-integrable perturbations to {H}ilbert--{S}chmidt operators},
  journal = {SIAM Journal on Mathematical Analysis},
  volume  = {57},
  number  = {6},
  pages   = {6103--6118},
  year    = {2025},
  doi     = {10.1137/24M1703136}
}

@article{Borcea2002,
  author  = {Borcea, Liliana},
  title   = {Electrical impedance tomography},
  journal = {Inverse Problems},
  volume  = {18},
  number  = {6},
  pages   = {R99--R136},
  year    = {2002},
  doi     = {10.1088/0266-5611/18/6/201}
}

@article{CaroRogers2016,
  author  = {Caro, Pedro and Rogers, Keith M.},
  title   = {Global uniqueness for the {C}alder{\'o}n problem with {L}ipschitz conductivities},
  journal = {Forum of Mathematics, Pi},
  volume  = {4},
  pages   = {e2},
  year    = {2016},
  doi     = {10.1017/fmp.2015.9}
}

@article{ClopFaracoRuiz2010,
  author  = {Clop, Albert and Faraco, Daniel and Ruiz, Alberto},
  title   = {Stability of {C}alder{\'o}n's inverse conductivity problem in the plane for discontinuous conductivities},
  journal = {Inverse Problems and Imaging},
  volume  = {4},
  number  = {1},
  pages   = {49--91},
  year    = {2010},
  doi     = {10.3934/ipi.2010.4.49}
}

@article{DaudeMaciaMeronoNicoleau2026,
  author  = {Daud{\'e}, Thierry and Maci{\`a}, Fabricio and Mero{\~n}o, Crist{\'o}bal and Nicoleau, Fran{\c{c}}ois},
  title   = {Stable factorization of the {C}alder{\'o}n problem via the {B}orn approximation},
  journal = {Journal of the London Mathematical Society (2)},
  volume  = {113},
  number  = {6},
  pages   = {e70575},
  year    = {2026},
  doi     = {10.1112/jlms.70575}
}

@article{DiCristoRondi2003,
  author  = {Di Cristo, Michele and Rondi, Luca},
  title   = {Examples of exponential instability for inverse inclusion and scattering problems},
  journal = {Inverse Problems},
  volume  = {19},
  number  = {3},
  pages   = {685--701},
  year    = {2003},
  doi     = {10.1088/0266-5611/19/3/313}
}

@article{Dobson1992,
  author  = {Dobson, David C.},
  title   = {Estimates on resolution and stabilization for the linearized inverse conductivity problem},
  journal = {Inverse Problems},
  volume  = {8},
  number  = {1},
  pages   = {71--81},
  year    = {1992},
  doi     = {10.1088/0266-5611/8/1/005}
}

@article{GardeHyvonen2024,
  author  = {Garde, Henrik and Hyv{\"o}nen, Nuutti},
  title   = {Linearized {C}alder{\'o}n problem: reconstruction and {L}ipschitz stability for infinite-dimensional spaces of unbounded perturbations},
  journal = {SIAM Journal on Mathematical Analysis},
  volume  = {56},
  number  = {3},
  pages   = {3588--3604},
  year    = {2024},
  doi     = {10.1137/23M1609270}
}

@article{GardeKnudsen2017,
  author  = {Garde, Henrik and Knudsen, Kim},
  title   = {Distinguishability revisited: depth dependent bounds on reconstruction quality in electrical impedance tomography},
  journal = {SIAM Journal on Applied Mathematics},
  volume  = {77},
  number  = {2},
  pages   = {697--720},
  year    = {2017},
  doi     = {10.1137/16M1072991}
}

@article{Isaacson1986,
  author  = {Isaacson, David},
  title   = {Distinguishability of conductivities by electric current computed tomography},
  journal = {IEEE Transactions on Medical Imaging},
  volume  = {5},
  number  = {2},
  pages   = {91--95},
  year    = {1986},
  doi     = {10.1109/TMI.1986.4307752}
}

@article{KnudsenLassasMuellerSiltanen2009,
  author  = {Knudsen, Kim and Lassas, Matti and Mueller, Jennifer L. and Siltanen, Samuli},
  title   = {Regularized {D}-bar method for the inverse conductivity problem},
  journal = {Inverse Problems and Imaging},
  volume  = {3},
  number  = {4},
  pages   = {599--624},
  year    = {2009},
  doi     = {10.3934/ipi.2009.3.599}
}

@misc{LiimatainenSaloTrad2026,
  author        = {Liimatainen, Tony and Salo, Mikko and Trad, William},
  title         = {Improved stability of low {F}ourier modes in inverse problems for potentials},
  year          = {2026},
  note          = {Preprint, arXiv:2608.20673},
  eprint        = {2608.20673},
  archiveprefix = {arXiv},
  primaryclass  = {math.AP}
}

@misc{Lionheart2026,
  author        = {Lionheart, William R. B.},
  title         = {Explicit diagonalization of the linearized {C}alder{\'o}n problem on the disk: multiplicity normalization and exponential spectral attenuation},
  year          = {2026},
  note          = {Preprint, arXiv:2609.23740},
  eprint        = {2609.23740},
  archiveprefix = {arXiv},
  primaryclass  = {math.AP}
}

@misc{LionheartPidcockSebu2026,
  author        = {Lionheart, William R. B. and Pidcock, Michael and Sebu, Cristiana},
  title         = {At last, an explicit singular value decomposition for electrical impedance tomography},
  year          = {2026},
  note          = {Preprint, arXiv:2610.35217},
  eprint        = {2610.35217},
  archiveprefix = {arXiv},
  primaryclass  = {physics.med-ph}
}

@article{Mandache2001,
  author  = {Mandache, Niculae},
  title   = {Exponential instability in an inverse problem for the {S}chr{\"o}dinger equation},
  journal = {Inverse Problems},
  volume  = {17},
  number  = {5},
  pages   = {1435--1444},
  year    = {2001},
  doi     = {10.1088/0266-5611/17/5/313}
}

@article{NagayasuUhlmannWang2009,
  author  = {Nagayasu, Sei and Uhlmann, Gunther and Wang, Jenn-Nan},
  title   = {A depth-dependent stability estimate in electrical impedance tomography},
  journal = {Inverse Problems},
  volume  = {25},
  number  = {7},
  pages   = {075001},
  year    = {2009},
  doi     = {10.1088/0266-5611/25/7/075001}
}

@article{KowWang2022,
  author  = {Kow, Pu-Zhao and Wang, Jenn-Nan},
  title   = {Refined instability estimates for some inverse problems},
  journal = {Inverse Problems and Imaging},
  volume  = {16},
  number  = {6},
  pages   = {1619--1642},
  year    = {2022},
  doi     = {10.3934/ipi.2022017}
}

@article{KowSaloZou2025,
  author  = {Kow, Pu-Zhao and Salo, Mikko and Zou, Sen},
  title   = {Increasing resolution and instability for linear inverse scattering problems},
  journal = {Journal of Functional Analysis},
  volume  = {289},
  number  = {1},
  pages   = {110923},
  year    = {2025},
  doi     = {10.1016/j.jfa.2025.110923}
}

@article{Rondi2006,
  author  = {Rondi, Luca},
  title   = {A remark on a paper by {A}lessandrini and {V}essella},
  journal = {Advances in Applied Mathematics},
  volume  = {36},
  number  = {1},
  pages   = {67--69},
  year    = {2006},
  doi     = {10.1016/j.aam.2004.12.003}
}

@article{SylvesterUhlmann1987,
  author  = {Sylvester, John and Uhlmann, Gunther},
  title   = {A global uniqueness theorem for an inverse boundary value problem},
  journal = {Annals of Mathematics (2)},
  volume  = {125},
  number  = {1},
  pages   = {153--169},
  year    = {1987},
  doi     = {10.2307/1971291}
}

@article{Uhlmann2009,
  author  = {Uhlmann, Gunther},
  title   = {Electrical impedance tomography and {C}alder{\'o}n's problem},
  journal = {Inverse Problems},
  volume  = {25},
  number  = {12},
  pages   = {123011},
  year    = {2009},
  doi     = {10.1088/0266-5611/25/12/123011}
}

@article{KochRulandSalo2025,
  author  = {Koch, Herbert and R{\"u}land, Angkana and Salo, Mikko},
  title   = {On instability mechanisms for inverse problems},
  journal = {Ars Inveniendi Analytica},
  year    = {2025},
  note    = {Revised Paper No.~1, 94 pp.},
  doi     = {10.15781/kvf0-ds86},
  url     = {https://arxiv.org/abs/2012.01855v4}
}

@misc{OpenAIConductivity2026,
  author = {{OpenAI}},
  title  = {Nonuniqueness for bounded measurable scalar conductivities
            in three dimensions},
  year   = {2026},
  note   = {Preprint, September 23.
            \href{https://github.com/openai/math/blob/main/preprints/Nonuniqueness-for-Bounded-Measurable-Scalar-Conductivities-in-Three-Dimensions-September-23-2026/paper.pdf}{Full text}}
}
\end{document}